\documentclass{article}
\usepackage{graphicx}
\usepackage{amsfonts,amsmath,amssymb,amsthm, amscd}

\input{epsf}

\numberwithin{equation}{section}
\theoremstyle{definition}

\newtheorem{theorem}{Theorem}[section]
\newtheorem{proposition}[theorem]{Proposition}
\newtheorem{lemma}[theorem]{Lemma}
\newtheorem{corollary}[theorem]{Corollary}

\newtheorem{definition}[theorem]{Definition}
\newtheorem{example}[theorem]{Example}
\newtheorem{remark}[theorem]{Remark}

\def\Z{\mathbb Z}
\def\C{{\mathbb C}}
\def\A{{\cal A}}
\def\ni{\noindent}

\def\G{\Gamma}

\def\<{\langle}
\def\>{\rangle}

\def\g{\gamma}
\def\d{\delta}

\def\t{t^{\pm 1}}

\def\D{\Delta}
\def \={ {\buildrel \cdot \over =}}

\begin{document}

\title{Crowell's Derived Group and Twisted  Polynomials}

\author{Daniel S. Silver \and Susan G. Williams\thanks{Both authors partially supported by NSF grant
DMS-0304971.} \\ {\em
{\small Department of Mathematics and Statistics, University of South Alabama}} \\ {\em
{\small silver@jaguar1.usouthal.edu, swilliam@jaguar1.usouthal.edu}}}

\maketitle %{\setlength{\linewidth}{2in}

\centerline{\it Dedicated to Louis H. Kauffman on the occasion of his sixtieth birthday.}

\begin{abstract} The derived group of a permutation representation, introduced by R.H. Crowell, unites many notions of knot theory. We survey Crowell's construction, and offer new applications. 

The twisted Alexander group of a knot is defined. Using it, we obtain twisted Alexander modules and polynomials. Also, we extend a well-known theorem of Neuwirth and Stallings giving necessary and sufficient conditions for a knot to be fibered.

Virtual Alexander polynomials provide obstructions for a virtual knot that must vanish if the knot has a diagram with an Alexander numbering. The extended group of a virtual knot is defined, and using it a more sensitive obstruction is obtained. 
\end{abstract} 

\begin{footnote}{Mathematics Subject Classification:  
Primary 57M25; secondary 20F05, 20F34.}\end {footnote}

\section{Survey of Crowell's derived group.} \label{survey} 

The derived group of a permutation representation unifies in a single theory the concepts of knot group, Alexander matrix and covering space, three concepts that Crowell and Fox believed were central to knot theory (see preface to \cite{CF77}).  Crowell presented his idea in \cite{C84} as a culmination of many years of thought. 

The operator notation that we use is more compact than the notation of \cite{C84}. It is convenient for our applications, and we hope that it will suggest further applications. 

Throughout $G$ will denote a multiplicative group acting on the right of a nonempty set $\G$. The action $\G\times G \to \G$ is denoted by $(\g,g)\mapsto \g g$. It determines in the usual way a representation $\rho: G \to S_\G$, where $S_\G$  is the group of permutations of $\G$. 

\begin{definition} \cite{C84} The \emph{derived group (of the permutation representation $\rho$)} is the free group with basis $\G\times G = \{ g^\g\mid g\in G, \g\in \G\}$ modulo the relations
$(gh)^\g = g^\g h^{\g g}$, for all $g, h \in G$ and $\g \in \G$. We denote it by $G_\rho$. \end{definition}

It is helpful to regard the exponent $\g$ in $g^\g$ as a coordinate.
Before giving examples, we mention basic combinatorial facts about derived groups that are proved in \cite{C84}. 

\begin{lemma}  \label{combinatorial} (i) $g^\g = 1$ if and only if $g$ is trivial in $G$.

\ni (ii) $(g^\g)^{-1} = (g^{-1})^{\g g}$, for any $g \in G$.

\ni (iii) Every nontrivial element $u \in G_\rho$ has a one and only one \emph{normal form} 
$$u= g_1^{\g_1}g_2^{\g_2}\cdots g_n^{\g_n},$$
where $g_k \ne 1$ and $\g_{k+1} \ne \g_k g_k$, for $1\le k < n$. \end{lemma}

As Crowell warned, $G_\rho$ should not be confused with the commutator subgroup $[G, G]$ (also known as the derived group of $G$). The reason for the name is a connection between $G_\rho$ and the derived module of a homomorphism.

\begin{definition} \cite{C71} \cite{C84} Let $P: G \to \G$ be a group homomorphism. Let $\G \times G \to \G$ be the $G$-action on $\G$ given by
$(\g, g) \mapsto \g P(g)$, and let  $\rho: G \to S_\G$ be the associated permutation representation.  The \emph{derived module} of $P$ is the $\Z[\G]$-module $G_\rho/ [G_\rho, G_\rho]$. \end{definition}

Derived groups can be regarded in two different ways as universal objects in appropriate categories. We briefly outline each approach.

Let $A$ be a multiplicative group. We recall that a \emph{crossed product} is a function $f: \G \times G \to A$ with the property that $f(\g, g_1g_2)= f(\g, g_1)f(\g g_1, g_2)$, for each $\g \in \G,\  g_1, g_2 \in G$.  The function $\wedge: \G \times G \to G_\rho$ determined by 
$(\gamma, g) \mapsto g^\g$ is an example of a crossed product. Moreover, given any crossed product $f: \G\times G \to A$ there is a unique group homomorphism $\bar f: G_\rho \to A$ such that $ \bar f \wedge =f.$  This universal property characterizes $G_\rho$ in the usual way. 

A second point of view, one that is more topological, involves covering spaces. Let $B, b_0$ be a connected space  and  base point such that $\pi_1(B, b_0) \cong G$. Let $p: E \to B$ be the covering space corresponding to the permutation representation $\rho: G \to S_\G$. The preimage $p^{-1}(b_0)$ can be identified with $\G$. Let $(\g, g) \in \G\times G$ denote the relative homotopy class of paths lying above $g$ beginning at $\g$ and ending at $\g g$. Then $\G \times G$ is a groupoid under concatenation, and 
$$(\g_1, g_1) (\g_1g_1, g_2) = (\g_1, g_1g_2),$$
for $g_1, g_2 \in G, \ \g \in \G$. The function $\wedge: \G \times G \to G_\rho$ sending $(\g, g)$ to $g^\g$ is a morphism of groupoids. It has the universal property that for any group $A$ and groupoid morphism $f: \G \times G \to A$, there is a unique homomorphism $\bar f: G_\rho \to A $ such that $\bar f \wedge = f$. 

In a sense, $G_\rho$ is the smallest group generated by the monoid of path liftings. It is shown in \cite{C84} that $G_\rho$ is isomorphic to the free product of $\pi_1(E)$ and a free group $F$. When $B$ is a cell complex, free generators of $F$ can be identified with the edges of a maximal tree in the $1$-skeleton of the induced cell complex of $E$. 

Under certain conditions the derived group has a natural structure as a group with operators. 
%Operator groups are both natural and useful when considering fundamental groups of covering spaces. 
The notion, reviewed here, was developed by Noether\cite{N29}, who attributed the idea to Krull.  A detailed treatment can be found in \cite{B74}. 
An \emph{ $\Omega$-group} is a group $K$ together with a set $\Omega$ (\emph{operator set}) and a function $\phi: K \times \Omega \to K,\ (g, \omega)\mapsto g^\omega$, such that for each fixed $\omega$, the restricted map 
$\phi_\omega: K \to K$ is an endomorphism. 
If $K$ is an $\Omega$-group and $\Omega_0$ is a subset of $\Omega$, then we can regard $K$ as an $\Omega_0$-group by restricting the action.

An $\Omega$-group is \emph{ finitely generated}
(respectively, \emph{ finitely presented}) if it is generated (resp. presented) by finitely many $\Omega$-orbits of generators (resp. generators and relators). The commutator subgroup of a knot group is an example of a finitely presented $\Z$-group. Free $\Omega$-groups are defined in the standard way.

Assume that a group $G$ acts on a multiplicative semigroup $\G$ so that 
$(\g \g') g = \g (\g' g)$ for all $\g, \g' \in \G$ and $g \in G$. Let $G_\rho$ be the associated derived group. We can define a map $\phi: G_\rho \times \G \to G_\rho$ on generators of $G_\rho$ by $\phi_\g(g^{\g'}) = g^{\g \g'}$, extending multiplicatively. One checks that $G_\rho$ is then a $\G$-group. 
Its abelianization is a left $\Z[\G]$-module with 
$(\sum n_i \g_i)g^\g = \sum n_i g^{\g_i\g}$ for $\sum n_i \g_i \in {\Bbb Z}[\Gamma]$. If $G$ has presentation $\<g_1, \ldots, g_m \mid r_1,\ldots, r_n\>$, then $G_\rho$ is finitely presented as a $\G$-group by the orbits of $g_1^1, \ldots, g_m^1$ and $r_1^1, \ldots, r_n^1$. We will ususally abbreviate $g_i^1$ by $g_i$.  The relators $r_j^1$ are easily written in terms of generators $g_i^\g$ of $G_\rho$, as we see in the following example. 

\begin{example} \label{figure8} The group $\pi_k$ of the figure-eight knot $k=4_1$ has presentation 
\begin{equation} \label{pi}  \pi_k=\< a, b \mid \bar a b a \bar b a b \bar a \bar b a \bar b\>,
\end{equation} 
where $\ \bar{}\ $ denotes  inverse. Let $P: \pi_k \to \<t\mid\ \>\cong \Z$ be  the abelianization homomorphism mapping  $a$ and $b$ to $t$.  The derived group of the associated permutation representation $\rho:\pi_k \to S_\Z$ has $\Z$-group presentation 
\begin{equation*} \pi_{k,\rho}=\< a, b \mid \bar a^{t^{- 1}} b^{t^{- 1}} a \bar b a b^t \bar a^t\bar b a \bar b\>.\end{equation*}
Here the $\Z$-group generators $a$ and $b$ represent orbits $\{a^{t^i}\}_{i\in \Z}$ and $\{b^{t^i}\}_{i\in \Z}$ of group generators, while  the single $\Z$-group relator $\d(r)$ represents an orbit $\{r^{t^i}\}_{i\in \Z}$ of group relators.

The relator $\d(r)$ is a rewrite of $r^1$, where $r$ is the relator in \ref{pi}. It is easily computed by a variant of Fox calculus, using the axioms:
$$\d(g) =g, \quad \d(\bar g)=\bar g^{P(\bar g)}$$
$$\d(wg) = \d(w)g^{P(w)},\quad \d(w\bar g) = \d(w)\bar g^{P(w\bar g)}.$$
Here $g$ is any generator while $w$ is a word in generators and their inverses.  It is easy to check that $r^{t^i}$ can be rewritten as $\delta(r)^{t^i}$, that is, $\delta(r)$ with each generator $g^{t^j}$ replaced by $g^{t^{i+j}}$.

For this particular representation $\rho$, we refer to the group $\pi_{k,\rho}$ as the \emph{ Alexander group} of the knot,
and denote it by $\A_k$.  Abelianizing $\A_k$ gives the derived module of the abelianization homomorphism. It is the Alexander module of the knot, a  $\Z[\t]$-module with matrix presentation described by the usual Fox partial derivatives (see \cite{BZ03}):

$$\Biggr( \biggr( \frac{\partial r} {\partial a}\biggr)^P\ \biggr(\frac{\partial r} {\partial b}\biggr)^P\Biggr) = (-t^{-1}+3 -t  \ \ \quad t^{-1}-3 +t).$$
The superscript indicates that each term of the formal sum is to be replaced by its $P$-image. 

The $0$th characteristic polynomial of the matrix, $t^2-3t +1$ (well defined up to multiplication by $\pm t^i$), is the ($1$st) Alexander polynomial $\D_k(t)$ of the knot. The index shift is a consequence of the fact that we removed a redundant relator from the Wirtinger presentation for $\pi_k$ when obtaining  \ref{pi}.

Replacing $\G = \Z$ by $\G = \Z/r\Z \cong \<t\mid t^r\>$ results in another derived group that is identical to the one above except that exponents $t^i$ are taken modulo $r$. In this case, $G_\rho$ modulo the normal subgroup generated by $x, x^t, \ldots, x^{t^{r-2}}$ produces the fundamental group of the $r$-fold cyclic cover of ${\Bbb S}^3 \setminus k$. Killing the entire family of generators,
$x, x^t, \ldots, x^{t^{r-1}}$ results in the fundamental group of the 
$r$-fold cover $M_r$ of ${\Bbb S}^3$ branched over $k$. 

\end{example}
\section{Core group as a derived group.}
\begin{figure}
\begin{center}
\includegraphics[width=.8in]{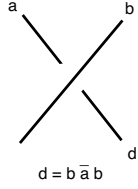}
\caption{Relations for core group $C_k$.}
\label{core}
\end{center}
\end{figure}

A presentation for the \emph{core group} $C_k$ of a knot can be obtained from any diagram. Generators correspond to arcs, while at each crossing the generators must satisfy a relation indicated in Figure \ref{core}. The relation does not depend on the orientation of the arcs. 

The core group was discovered independently by A.J. Kelly \cite{kelly} and M. Wada \cite{W92}. It appears in a paper by R.A. Fenn and C.P. Rourke \cite{FR92}. Wada proved in \cite{W92} that 
$C_k$ is isomorphic to the free product of $\pi_1 M_2$ with an infinite cyclic group (see \cite{P98} for an elementary proof).

The \emph{$\pi$-orbifold group} $O_k$ is another group associated to any knot $k$. It is the quotient group 
$$\pi_k/ \<\<x^2\>\>, $$
where $x$ is any meridian element and $\<\<\quad  \>\>$ denotes normal closure. Letting $\tau$ denote the nontrivial covering transformation of $M_2$, the $\pi$-orbifold group has a topological interpretation as the group of lifts of $\tau$ and the identity transformation to the universal cover of $M_2$. Consequently, it fits into a short exact sequence
$$1 \to \pi_1 M_2 \to O_k\ {\buildrel P \over \to}\ \Z/2 \to 1$$
(see for example \cite{kawauchi}, page 133). The homomorphism $P$ induces a permutation representation $\rho: O_k \to S_{\Z/2}$. Hence the derived group $O_{k,\rho}$ is defined.

\begin{proposition} The derived group $O_{k,\rho}$ is isomorphic to the core group $C_k$. 

\end{proposition}

\begin{proof} Consider a Wirtinger presentation for the group of $k$,  $$\pi = \< a, b, \ldots \mid r, s, \ldots \>,$$
where as usual generators correspond to arcs of a diagram, and relations correspond to crossings (see Figure \ref{wirtinger}).
Since any two (meridian) generators are conjugate, the $\pi$-orbifold group has presentation 
$$O_k = \< a, b, \ldots \mid r, s, \ldots, a^2, b^2, \ldots \>.$$
The homomorphism $P: O_k \to \Z \cong \<t\mid t^2\>$ takes each 
generator to $t$. Hence the derived group of the corresponding permutation representation $\rho$ has group presentation 
$$O_{k,\rho} = \< a^{t^i}, b^{t^i}, \ldots \mid \d(r)^{t^i}, \d(s)^{t^i},  \ldots, \d(a^2)^{t^i},  \d(b^2)^{t^i}, \ldots\>,$$
where $i$ is taken modulo $2$, and as before $\d(r), \d(s), \ldots$ indicate rewritten relators.

Consider a typical Wirtinger relation $bd = ab$ corresponding to a positve crossing, as in Figure \ref{wirtinger}. (The argument or a negative crossing is similar.) The relation determines a pair of relations in the presentation of the derived group. Abbreviating $a^1, b^1, \ldots$ by $a, b, \ldots,$ the relations appear as:
\begin{equation} \label{core1}  b d^t = a b^t, \ b^t d = a^t b.
\end{equation} 

The relators $a^2, b^2, \ldots$ also determine pairs of relators:
$$a a^t, \ a^t a, \ b b^t = b^t b, \ldots$$
or equivalently, 
\begin{equation} \label{core2} a^t = \bar a, \ b^t = b, \ldots.
\end{equation}
Using the relations \ref{core2}, we can eliminate $a^t, b^t, \ldots,$ from the presentation for $O_{k, \rho}$. The relations \ref{wirtinger} become 
$b \bar d = a \bar b$ (the second relation being redundant).
Since this relation is equivalent to $d = b \bar a  b$, we have arrived at a presentation of the core group. Hence $O_{k, \rho} \cong C_k$.

\begin{figure}
\begin{center}
\includegraphics[width=1.8in]{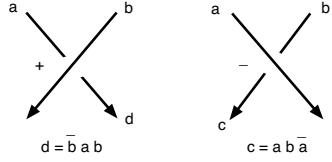}
\caption{Wirtinger relations.}
\label{wirtinger}
\end{center}
\end{figure}

\end{proof}

\section{Twisted invariants.}\label{twisted}

Suppose the permutation representation $\rho$ of $G$ corresponds to an action of $G$ on a left $R$-module $V$. Then, regarding $R$ as a multiplicative semigroup, the derived group $G_\rho$ has the structure of an $R$-group. The requisite function $G_\rho \times R \to G_\rho$ is defined on generators by $(g^v, r) \mapsto  g^{r v}$ and extended multiplicatively.  This is well defined on $G_\rho$ since $r$ takes $(gh)^v=g^vh^{v\rho(g)}$ to $(gh)^{rv}=g^{rv}h^{rv\rho(g)}$.

Note that in $G_\rho$ we do \emph{not} have the relation  $g^vg^w=g^{v+w}$. If we impose this relation for all $g\in G$ and $v, w\in V$ we obtain a quotient group with an induced $R$-group structure.  Its abelianization is an $R$-module, since we recover the additivity of the $R$ action from
\begin{equation*}\begin{split}(r_1+r_2)(vg) & =((r_1+r_2)v)g=(r_1v+r_2v)g\\  &=(r_1v)g+(r_2v)g=r_1(vg)+r_2(vg),\end{split}\end{equation*}
where the third equality follows from the imposed relations.

\begin{example}\label{twistedfigure8} 
Twisted Alexander polynomials, introduced for knots in \cite{L90} and defined for finitely presented groups in \cite{W94}, can be obtained by replacing $\rho$ with more general representations.
We return to the example above of the figure-eight knot (example \ref{figure8}). 
Let $P_1: \pi_k \to \<t\mid \ \>$ be the abelianization homomorphism considered there, and let $P_2: \pi_k \to {\rm SL}_2(\C)$ be the discrete faithful representation 
$$a \mapsto \begin{pmatrix}1\ &1\ \\ 0\ &1\ \end{pmatrix},\quad  b\mapsto \begin{pmatrix}1\ &0\  \\ -w\ &1\ \end{pmatrix},$$
where $\omega=\displaystyle{ \tfrac{1-\sqrt{-3}} 2}$ (see Chapter 6 of \cite{T78}). 
Together, $P_1$ and $P_2$ determine a homomorphism $P:\pi_k \to 
GL_2(\C[t^{\pm 1}])$ described  by:
$$a \mapsto  \begin{pmatrix} t\ &t\ \\ 0\  &t\  \end{pmatrix} \quad b\mapsto \begin{pmatrix}t&0\\-t\omega&t \end{pmatrix}. $$
%One can regard $P(g)$ as $P_1(g)\otimes P_2(g)$, the product of the scalar monomial $P_1(g)$ with $P_2(g)$ in $\Z[\t]\otimes_\Z{\rm SL}_2(\C) \subset {\rm GL}_2(\C[\t])$. 
Denote the image of $P$ by $\G$.

The derived module of $P:\pi_k \to \G$ is a $\Z[\G]$-module with presentation matrix given by a $1 \times 2$ Jacobian: 
$$\Biggr( \biggr(\frac{\partial r}{\partial a}\biggr)^P \ \biggr(\frac{\partial r}{\partial b}\biggr)^P\Biggr).$$
Each entry is a formal sum of matrices. Add its terms to form a single $2\times 2$-matrix. The result is: 
\begin{equation} \Biggr( \begin{pmatrix} -\bar \omega t +1 +\bar \omega -\tfrac{1}{ t} &\quad  t-2-\omega +\tfrac{1}{ t}\\  i \sqrt{3} \bar \omega (t-1) & -\bar \omega t + 3 -\tfrac{1}{ t}\end{pmatrix} \quad \begin{pmatrix} -2+\tfrac{1}{ t} &\quad  -\omega t + 1 +\omega -\tfrac{1}{ t}\\ \bar \omega t -1 &\quad 2t-3+\tfrac{1}{ t}\end{pmatrix}\Biggr).\end{equation} \label{2by2}
Wada \cite{W94} showed that the following procedure produces a polynomial invariant of knots.  Remove the inner parentheses to form  a $2\times 4$ \emph{ twisted Alexander matrix} $A_\rho$.
Delete the first two columns (corresponding to the generator $a$), take the determinant, and divide by ${\rm Det}(tI - P(a))$. The resulting polynomial is Wada's invariant,
$$W_{k, \rho} = \frac  {(t-1)^2 (t^2 - 4 t +1)} {(t-1)^2} = t^2 - 4t +1.$$
The result is the same if we reverse the roles of $a$ by $b$. (See \cite{W94}.) 

Adding terms and removing inner parentheses to obtain $B$ appears {\it ad hoc}. The following alternative approach, based on the derived group of an action on an $R$-module, gives $B$ more naturally.  

The group $\pi_k$ acts on the free $\C[\t]$-module $V = \C[\t]\oplus \C[\t]$ 
via $P: \pi_k \to {\rm GL}_2(\C[\t])$: for any $v \in V$ and $g \in \pi_k$, we define $vg$ to be $vP(g)$. This gives a permutation representation $\rho_V$ of $\pi_k$ on $V$. By remarks above, the derived group $\pi_{\pi_k, \rho_V}$ is a $\C[\t]$-group, and the abelianization becomes a left $\C[\t]$-module if we impose the extra relations $(v_1 + v_2) g = v_1 g + v_2 g$ for every $v_1, v_2 \in V$ and $g\in \pi_k$. It is generated by $e_1 a, e_2 a, e_1 b, e_2 b$, where $e_1=(1,0), e_2=(0,1)$ is the standard basis for $V$, and one easily checks that $A_\rho$ is a presentation matrix for it. \end{example}

\begin{remark} (i) Twisted Alexander polynomials of any $d$-component link $\ell$ can be found in a similar way, using the abelianization $P_1: \pi_\ell \to \Z^d$. 

(ii)  The greatest common divisor of the six $2\times 2$ minors of $A$ is an invariant of $k$ and the representation $P$.  It is the $1$st twisted Alexander polynomial $\D_1$, defined by P. Kirk and C. Livingston \cite{KL99}. 
The greatest common divisor of the $1\times 1$ minors (that is, the entries of $A$) is the $0$th twisted polynomial  $\D_0$. The relationship
$W_{k, \rho} = \D_1/\D_0$ holds generally. \end{remark}

The reader is cautioned that an invariant does not result by simply taking the greatest common divisor of the two $2\times 2$ minors corresponding to the entries of \ref{2by2}, contrary to Proposition 1.3 of \cite{JW93}. Such a quantity is not preserved when the presentation of $\pi_k$ is altered by a Tietze move, in contradiction to the assertion on page 213 of \cite{JW93}. Indeed, in the above example, the greatest common divisor of the two minors is $(t-1)^2(t^2-4t-1)$. However, if in the presentation of $\pi_k$ we replace $a$ by the new generator $c= \bar a b$, accomplished by a sequence of Tietze moves, then the polynomial $t^2-4t+1$ is obtained.

The Alexander group $\A_k$ is the free product of the commutator subgroup $[\pi_k, \pi_k]$ with a countable-rank free group $F$. One can choose $F$ to be the free group generated by any $\Z$-group family $\{x^{t^i}\}$ of group generators, where $x$ is a Wirtinger generator (cf. Section \ref{survey}; also see \cite{SW01}). Notice, in particular, that $\A_k/\< \< x \> \>$  is isomorphic to $[\pi_k, \pi_k]$. 
Here $\<\< x \>\>$ denotes the normal  $\Z$-subgroup generated by $x$, that is, the smallest normal subgroup containing the orbit $\{x^{t^i}\}$. 

For any knot $k$, the commutator subgroup $[\pi_k, \pi_k]$ is a finitely generated $\Z$-group. As an ordinary group, it need not be finitely generated.
A theorem of J. Stallings \cite{S62} states that if it is, then $k$ is fibered, and $[\pi_k, \pi_k]$ is free of rank $2 g_k$, where $g_k$ denotes the genus of $k$. 
The converse is also true: if $k$ is fibered, then $[\pi_k, \pi_k]$ is a free group of rank $2g_k$  \cite{N63}. 

A consequence is that if $k$ is fibered, then its Alexander polynomial is monic with degree equal to twice the genus of $k$.  H. Goda, T. Kitano and T. Morifuji \cite{GTM03} extended this classical result by showing that if $k$ is fibered and $\rho$ is obtained as in Example \ref{twistedfigure8} from the abelianization homomorphism $P_1$ and homomorphism $P_2: \pi_k \to SL_{2n}({\Bbb F}),\ {\Bbb F}$ a field, then Wada's invariant $W_{k,\rho}(t)$ is a rational function of monic polynomials. A similar result was obtained by J.C. Cha \cite{C03}. We prove a companion result for Alexander groups, generalizing the theorem of Neuwirth and Stallings. 

Consider a knot $k\subset {\mathbb S}^3$ with abelianization homomorphism $P_1: \pi_k \to \<t\mid\ \>$ and $P_2: \pi_k \to
{\rm SL}_n(R)$ any homomorphism. We assume that $R$ is a unique factorization domain.
Let $P: \pi_k \to 
{\rm GL}_n(R[t^{\pm 1}])$ be the product homomorphism, as above, with image $\G$. Let $\rho: \pi_k \to S_\G$ be the associated permutation representation. We denote  by  ${\cal R}(\pi)$ the collection of  permutation representations that arise this way. 

The derived group $\pi_{k, \rho}$ is the \emph{ $\rho$-twisted Alexander group} of $k$, and we denote it by $\A_{k,\rho}$. When $P_2$ is trivial, $\pi_{k, \rho}$ reduces to the Alexander group of $k$. The twisted Alexander group is a $\G$-group. The image $\G_2$ of $P_2$ is a subgroup of $\G$. We can also regard $\A_{k, \rho}$ as a $\G_2$-group,  letting only elements of $\G_2$ act on it. Throughout the remainder of the section,  $x$ will denote a Wirtinger generator of $\pi_k$ corresponding to a meridian of $k$.

\begin{theorem}\label{fibered}  If a  knot $k$ is fibered, then, for any representation $\rho\in {\cal R}(\pi)$, the twisted Alexander group $\A_{k,\rho}$ modulo $\<\<x\>\>$ is a free $\G_2$-group of rank $2g_k$.  Conversely, if $\A_{k,\rho}$ modulo
$\<\<x\>\>$ is a finitely generated $\G_2$-group for some $\rho\in {\cal R}(\pi)$, then $k$  is fibered. \end{theorem}

\begin{proof}  Assume that $k$ is fibered. Then $\pi_k$ is an HNN extension of a free group $F=F(a_1, \ldots, a_{2g})$, where $g=g_k$ is the genus of $k$, and it has a presentation of the form 
$$\<x, a_1, \ldots, a_{2g}\mid x a_1 x^{-1} =\phi (a_1), \ldots, x a_{2g}x^{-1}=\phi (a_{2g})\>\eqno(2.3)$$
for some automorphism $\phi$ of $F$ (see \cite{BZ03}, for example). The twisted Alexander group $\A_{k,\rho}$ has  a $\Gamma$-group presentation
$$\< x, a_1, \ldots, a_{2g}\mid xa_1^{tP_2(x)} \bar x^{P_2(xa_1 \bar x)} = \widetilde {\phi(a_1)}, \ldots,\ xa_1^{tP_2(x)} \bar x^{P_2(xa_{2g} \bar x)} = \widetilde {\phi(a_{2g})}\>, $$ 
where $\ \widetilde{}\ $ denotes the $\d$-rewrite,  explained above. The quotient $\G$-group $\A_{k,\rho}/\<\<x\>\>$ has presentation 
$$\<a_1, \ldots, a_{2g}\mid a_1^{tP_2(x)}= \widetilde  {\phi(a_1)}, \ldots,\ a_{2g}^{tP_2(x)} = \widetilde  {\phi(a_{2g})}\>.$$
Since $P_2(x)$ is invertible, each $a_i^t$ can be expressed as
$\widetilde {\phi(a_i)}^{{P_2(x)}^{-1}}$. Generally, the relators can be used one at a time to express each $a_i^{t^j},\ j \ne 0$, as a word in the $\G_2$-orbits of $a_1,\ldots, a_{2g}$. Hence $\A_{k,\rho}/\<\<x\>\>$ is isomorphic to the free $\G_2$-group generated by $a_1, \ldots, a_{2g}$.

Conversely, assume that $\A_{k,\rho}/\<\< x \>\>$ is a finitely generated $\G_2$-group for some $\rho \in {\cal R}(\pi)$. The knot group $\pi_k$ is generated by $x, a_1, \ldots, a_m$, where $a_1, \ldots, a_m \in [\pi_k, \pi_k]$. Then $\A_{k,\rho}/\<\< x \>\>$ is generated by elements of the form
$a_i^{t^j M}$, where $j \in \Z$ and $M$ ranges over matrices in $\G_2$, while $[\pi_k, \pi_k]$ is generated by the elements $a_i^{t^j}$. The mapping
$a_i^{t^j M}\mapsto a_i^{t^j}$ determines a surjection 
$f: \A_{k,\rho}/\<\< x \>\>\to [\pi_k, \pi_k]$. By assumption, we can find generators $a_i^{t^j M}$ of  $\A_{k,\rho}/\<\< x \>\>$ such that the values of $j$ are bounded. Consequently, the image of $f$ is finitely generated. Hence $k$ is fibered. \end{proof}

%%%%%%%%%%%%%%%%%%%%%%%%%%%%%%%%%

\section{Virtual knots and links.}  Like their classical counterparts, virtual knots and  links are equivalence classes of diagrams, the equivalence relation generated by Reidemeister moves. However, in the case of virtual knots and links,  we also allow virtual crossings (indicated by a small circle about the crossing), and the set of Reidemeister moves is suitably generalized. The notion is due to Kauffman, and the reader is referred to \cite{K98} or \cite{K99} for details. A result of M. Goussarov, M. Polyak and O. Viro \cite{GPV00} assures us that if two classical diagrams---that is, diagrams without virtual crossings---are equivalent by generalized Reidemeister moves, then they are equivalent by the usual classical Reidemeister moves. In this sense, virtual knot theory extends the classical theory.

Many invariants of classical knots and  links are also defined in the  virtual category. 
The knot group $\pi_\ell$ is such an invariant. Given any diagram of a virtual knot or link $\ell$,
a Wirtinger presentation is obtained by assigning a generator to each arc. (An arc is a maximal connected component of the  diagram containing no classical  under-crossing.) A relation is associated to each classical crossing in the  usual way. 

Let $\ell=\ell_1\cup \cdots \cup \ell_d$ be an oriented virtual link of $d$ components with group $\pi_\ell$. Regard $\Z^d$ as a multipicative group freely generated by $u_1,\ldots, u_d$. Let $P: \pi_\ell \to \Z^d$ be the abelianization homomorphism mapping the class of the $i$th oriented meridian to $u_i$. Let $\rho$ be the associated permutation representation 
$\rho: \pi_\ell \to S_{\Z^d}$. The derived group is the Alexander group $\A_\ell$ introduced in \cite{SW01}. 

\begin{figure}
\begin{center}
\includegraphics[width=2in]{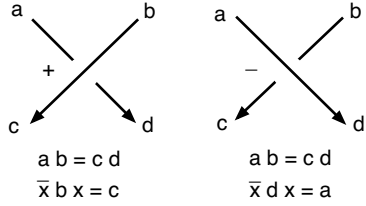}
\caption{Relations for extended group $\tilde \pi_\ell$.}
\label{fig1}
\end{center}
\end{figure}

Let $D$ be a diagram of $\ell$.  An \emph{edge} of $D$ is a maximal segment of an arc going from one classical crossing to the next. (Thus at each crossing, the arc passing over is broken into distinct edges.) The \emph{ extended group} $\tilde \pi_\ell$ has generators $a,b,c,\ldots$ corresponding to edges together with an additional generator $x$ not associated with any edge. Relations come in pairs, corresponding to classical crossings: $ab=cd,\ \bar xbx=c$, if the crossing is positive, and $ab =cd,\ \bar xdx=a$, if the crossing is negative (see Figure \ref{fig1}). It is a straightforward matter to check that the group  so defined is unchanged if a generalized Reidemeister move is applied to the diagram, and we leave this to the  reader.

The extended Alexander group $\tilde \A_\ell$ of the link $\ell$ was defined in \cite{SW01}
to be the $\Z^{d+1}$-group with generators $a,b,c,\ldots$, corresponding to edges as above, and relations at each classical crossing: $ab^{u_i} = c d^{u_j},\ c^v=b$, if the crossing is positive; $ab^{u_i} = c d^{u_j},\ a^v=d$, if the crossing is negative. We assume that  generators $a,d$ correspond to edges on the  $i$th component of the link while $b,c$ correspond to edges on the $j$th component. One regards each of $a,b,c,\ldots$ as families of generators indexed by monomials in $u_1, \ldots, u_d, v$,  generators for the group $\Z^{d+1}$ written multiplicatively. Similarly, each relation $r$ of the $\Z^{d+1}$-group presentation is a family 
of group relations indexed by monomials.

Let $P: \tilde \pi_\ell \to \Z^{d+1}$ be the abelianization homomorphism mapping generators corresponding to edges on the $i$th component to $u_i$, and mapping $x$ to $v$, and let $\rho: \tilde \pi_\ell \to S_{\Z^{d+1}}$ be the associated permutation representation. Consider the derived group $\tilde \pi_{\ell, \rho}$. It is generated as a $\Z^{d+1}$-group by $a, b,c,\ldots $ and $x$. Modulo $x$, the $\d$-rewrite of the relations for $\tilde \pi_{\ell, \rho}$ are easily seen to be those of the extended Alexander group. Indeed, the $\d$-rewrite of $ab=cd$ is 
$ab^{u_i} = cd^{u_j}$, while the rewrite of $xc\bar x = b$ is $xc^v  \bar x^{u_j} = b$. Killing $x$, that is, killing its $\Z^{d+1}$-orbit, yields the relations for $\tilde \A_\ell$ at a positive classical crossing. A similar argument applies to negative crossings. Summarizing:

\begin{proposition} Let $\ell=\ell_1\cup \cdots \cup \ell_d$ be an oriented virtual link of $d$ components, and let $\tilde \pi_\ell$ be its  extended group. Let $\tilde \pi_{\ell, \rho}$ be  the derived group of  the permutation representation associated to the abelianization homomorphism $P: \tilde \pi_\ell \to \Z^{d+1}$. Then as $\Z^{d+1}$-groups, the extended Alexander group $\tilde \A_\ell$ is isomorphic to $
\tilde \pi_{\ell, \rho}/\<\<x\>\>$, the derived group modulo the normal $\Z^{d+1}$-subgroup generated by $x$. \end{proposition}

Proposition \ref{invariant} yields new invariants for virtual links. 

\begin{proposition} \label{invariant} (i) Let $\ell$ be an oriented virtual link. Then
$$1\to \< \< x \> \> \to \tilde \pi_\ell\  {\buildrel p \over \to}\ \pi_\ell\to 1$$
is an exact sequence of groups, where $\<\<\ \  \>\>$ denotes normal closure and $p$ is the quotient map sending $x \mapsto 1$.\medskip

\ni (ii) If $\ell$ is a classical link, then
$$1\to \A_\ell \to \tilde \pi_\ell \ {\buildrel \chi \over \longrightarrow}\ \<t\mid\  \>\to 1$$
is an exact sequence of groups, where $\chi$ is the homomorphism mapping 
$x\mapsto t$  and $a,b,c, \ldots \mapsto1$.\end{proposition}

\begin{proof} Killing $x$ converts the relations for $\tilde \pi_\ell$ into the Wirtinger relations for $\pi_\ell$. Hence (i) is  proved. 

In order to prove (ii), assume that $D$ is a classical diagram for $\ell$; that is, a diagram without virtual crossings. By the Seifert smoothing algorithm, we can label the edges of $D$ with integers $\nu(a), \nu(b), \nu(c), \ldots$ such that at any crossing such as in Figure 1, 
whether positive or negative, edges corresponding to $a,c$ receive the same label, say $\nu \in \Z$, while edges corresponding to $b,d$ receive $\nu+1$. Such an assignment of integers will be called an \emph{ Alexander numbering}. (See \cite{SW01} for details.) 

After replacing each generator $a,b,c,\ldots$ with $x^{\nu(a)}ax^{-\nu(a)}, $ $x^{\nu(b)}bx^{-\nu(b)},$ $ x^{\nu(c)}cx^{-\nu(c)}\ldots$, the relations at a positive crossing become $a xb\bar x = c xd\bar x$ and $c=b$. At a negative crossing, the relations become
$a xb\bar x = c xd\bar x$ and $a=d$. Notice in particular that $\tilde \pi_\ell$ is generated by $x$ and symbols corresponding to the {\sl arcs} rather than the edges of $D$. Moreover, the form of the relations allows us to describe  the kernel of  $\chi$  using the  Reidemeister-Schreier method: it is generated as  a $\<t\mid\ \>$-group by symbols corresponding to the arcs of $D$ together with relations
$a b^t = b  d^t$ at a positive crossing and $ab^t=ca^t$ at a negative crossing. 
This is also a description of the Alexander group $\A_\ell$ of the  link. \end{proof}

Proposition \ref{invariant} motivates the following. 

\begin{definition} A virtual link is \emph{ almost classical} if some diagram for the link admits an Alexander numbering. \end{definition}

\begin{remark} (i) Reassuringly, classical implies almost classical. However, we see in Example \ref{almostclassical} that the converse does not hold. \medskip

(ii) The conclusion of Proposition \ref{invariant}(ii) holds for  almost classical links. The proof is similar. \end{remark}

\begin{figure}
\begin{center}
\includegraphics[width=1.7in]{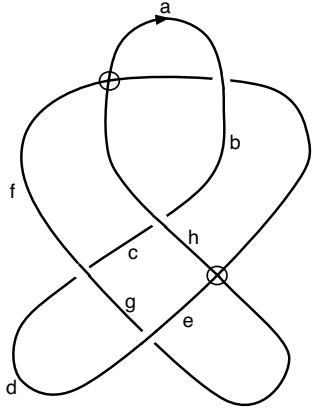}
\caption{Virtual knot $k$.}
\label{virtual5.2}
\end{center}
\end{figure}

\begin{example} \label{kishino} The virtual knot $k$ in Figure \ref{virtual5.2} appears in \cite{K02}.  Its extended group has generators $a,\ldots, h$. Using the relations corresponding to the classical crossings, one finds after some simplification that $\tilde \pi_k$ has group presentation
$$\<x,a,d\mid a x d \bar a \bar x \bar a x a^2 x a= d x a x d,\quad d x a x d = x a d x a\>.$$ The kernel of the homomorphism $\chi: \tilde \pi_k \to \<t\mid\ \>$ mapping $a,d \mapsto 1$ and $x \mapsto t$  is a $\<t\mid\  \>\cong \Z$-group. The following $\Z$-group presentation can be found using the Reidemeister-Schreier method: 
\begin{equation} \label{ker} {\rm ker}(\chi)= \<a, d\mid a d^t \bar a^t \bar a  (a^2)^t a^{t^2} = d a^t d^{t^2},\  d a^t d^{t^2} = a^t d^t a^{t^2}\> \end{equation}

If  $k$ were almost classical, then the group presented by \ref{ker} would be isomorphic to the Alexander group $\A_k$ of $k$. We compute a presentation of $\A_k$ directly from Figure \ref{virtual5.2}, and find
\begin{equation} \label{alex} \A_k = \<c,d \mid c d^t c^{t^2} = dc^t d^{t^2},\  c d^t c^{t^2} = dc^t d^{t^2}\>\end{equation} 
Each presentation (\ref{ker}), (\ref{alex}) determines a $\Z[\t]$-module with $2\times 2$-relation matrix. Both matrices have trivial determinant, which is the $0$th characteristic polynomial. However, the $1$st  characteristic polynomial corresponding to (\ref{ker}) is $1$ while that corresponding to (\ref{alex}) is $t^2-t+1$. Since these are unequal, the two groups are not isomorphic, and hence $k$ is not almost classical. \end{example}

\begin{remark} It follows that the knot in Example \ref{kishino} is not classical. Heather Dye has shown us an alternative proof using minimal surface representations, as in \cite{DK05}.

\end{remark}

The proof of Theorem 4.1 of  \cite{SW01} shows that if $\ell$ is an almost classical oriented virtual link of $d$ components, then a presentation for the Alexander group $\A_\ell$ obtained from any diagram of the link can be converted into a presentation for the extended Alexander group $\tilde \A_\ell$ by replacing each occurrence of $u_i$ by $u_i v,\ 1 \le i\le d$.  As a consequence, each virtual Alexander polynomial $\D_i(u_1, \ldots, u_d, v)$ is a polynomial in the $d$ variables $u_1v, \ldots, u_d v$. 

The arguments are similar for twisted groups and virtual Alexander polynomials.  Let $P_1: \tilde \pi_\ell \to \Z^{d+1}$ be the abelianization representation, and let 
$P_2: \tilde \pi_\ell \to {\rm GL}_n(\C)$ be a linear representation. Together, $P_1$ and  $P_2$ determine a representation 
$P: \tilde \pi_\ell \to {\rm GL}_n(\C[u_1^{\pm 1}, \ldots, u_d^{\pm 1}, v^{\pm 1}])$. As above, we denote the image of $P$ by $\G$ and the induced permutation representation of $\pi_\ell$ on $\G$ by $\rho$. A twisted Alexander matrix $A_\rho$ and associated polynomials $\D_{\ell, \rho, i}=\D_{\ell, \rho, i}(u_1, \ldots, u_d, v)$ can now be defined as in Example \ref{twistedfigure8}.

If the diagram for $\ell$ has $N$ classical crossings, $A_\rho$ is an $Nn\times Nn$-matrix. The greatest common divisor of the determinants of all $(Nn-i) \times (Nn-i)$ minors produces the  \emph{ $i$th twisted virtual Alexander polynomial} $\D_{\ell, \rho, i}$. As is the case with $\D_{\ell, 0}$, the twisted polynomial $\D_{\ell, \rho, 0}$ must vanish if $\ell$ is classical, and the argument is similar (see \cite{SW01}).

\begin{theorem} If $\ell$ is almost classical, then a presentation for the twisted Alexander group $\A_{\ell, \rho}$ obtained from any diagram of $\ell$ can be converted 
into a presentation for the extended Alexander group $\tilde \A_{\ell, \rho}$ by replacing each occurrence of $u_i$ by $u_iv,\ 1 \le i\le d$. \end{theorem}

\begin{corollary} Assume that $\ell$ is almost classical. Then each twisted virtual Alexander polynomial $\D_{\ell, \rho, i}(u_1, \ldots, u_d, v)$ is a polynomial in the $d$ variables $u_1v, \ldots, u_d v$. \end{corollary}

\section{Remarks about almost classicality.}
There are different ways to define what is means for a virtual knot to be almost classical. Our definition is motivated by group-theoretical concerns. The following two examples are intended to convince the reader that the definition is subtle and worthy of further study. 

\begin{figure}
\begin{center}
\includegraphics[width=2in]{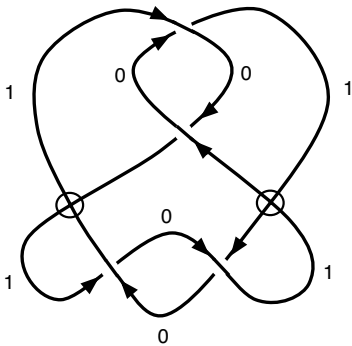}
\caption{ Non-classical knot that is almost classical.}
\label{ac5.6}
\end{center}
\end{figure}

\begin{example} \label{almostclassical} \emph{Almost classical does not imply classical.} Consider the diagram $D$ for a virtual knot $k$ in Figure \ref{ac5.6}.
An Alexander numbering for $D$ is shown, and hence $k$ is almost classical. The $1$st virtual Alexander polynomial $\D_1(u,v)$ is easily computed and seen to be $2-uv$. If $k$ were classical, then its Alexander polynomial would be $2-t$. 
However, $2-t$ is not reciprocal, and hence it cannot be the Alexander polynomial of any classical knot. From this we conclude that $k$ is not classical. \end{example}

One can weaken the definition of almost classical by asking only for a mod-$2$ Alexander numbering, that is, replacing $\Z$ by $\Z/2$ in the definition of Alexander numbering.   

\begin{figure}
\begin{center}
\includegraphics[width=1.5in]{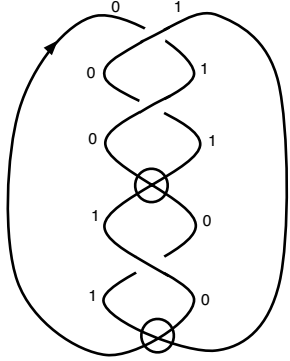}
\caption{Trefoil with virtualized crossing.}
\label{vtrefoil}
\end{center}
\end{figure}

\begin{example} \label{mod2} \emph{Existence of a mod-2 Alexander numbering does not imply almost classical.}  Figure \ref{vtrefoil} is a diagram of a virtual knot $k$ with a mod-$2$ Alexander numbering. The 0th virtual Alexander polynomial of $k$ is 
$$\D_{k, 0}(u,v) = 1 - u^2 -uv -v^2 +u^3v + uv^3 +u^2v^2 -u^3v^3.$$
Since it is not a polynomial in $uv$, the knot $k$ admits no diagram with an Alexander numbering. 

\end{example}

The knot in Example \ref{mod2} was introduced by Kauffman in 
\cite{K99}. It is obtained from a standard diagram of a trefoil by ``virtualizing a single crossing,'' switching the crossing and flanking it by virtual crossings. As observed in \cite{K99}, it is a nontrivial virtual knot with unit Jones polynomial. Such knots were further studied in \cite{SW04}.

%%%%%%%%%%%%%%%%%%%%%%%%%%%%%%%%%%%%%%

\end{document}